\documentclass{amsart}
\newcommand{\abs}[2][]{\mathopen#1|#2\mathclose#1|}
\newcommand{\norm}[2][]{\mathopen#1\|#2\mathclose#1\|}
\let\epsilon=\varepsilon
\let\phi=\varphi
\newtheorem{theorem}{Theorem}
\newtheorem{lemma}[theorem]{Lemma}
\let\Re\relax\DeclareMathOperator{\Re}{Re}

\begin{document}
\title{On the uniform convexity of \(L^p\)}
\author{Harald Hanche-Olsen}
\address{Norwegian University of Science and Technology\\
         Department of Mathematical Sciences\\
         NO--7491 Trondheim, Norway}
\email{hanche@math.ntnu.no}

\begin{abstract}
We present a short, direct proof of the uniform convexity of $L^p$
spaces for $1<p<\infty$.
\end{abstract}
\maketitle
\noindent
The standard proof of the uniform convexity of $L^p$ using Clarkson's
\cite{clarkson} or Hanner's \cite{hanner} inequalities
(see also \cite{balletal})
is rarely taught in functional analysis classes,
in part (the author imagines) because the proofs of those inequalities
are quite non-intuitive and unwieldy.
We present here a direct proof,
cheerfully sacrificing the optimal bounds --
for which, see \cite{hanner,balletal}.

In order to motivate our proof, and to explain how it might fit into a
standard functional analysis course,
we remind the reader of Young's inequality
\begin{equation} \label{young}
  \Re uv\le\frac{\abs{u}^p}{p}+\frac{\abs{v}^q}{q},
  \qquad 1<p<\infty,\ \frac1p+\frac1q=1.
\end{equation}
This has many proofs.
For example, if $u,v>0$ we write $u=e^{s/p}$ and $v=e^{t/q}$
and use the convexity of the exponential function to arrive at
$e^{s/p+t/q}\le e^s/p+e^t/q$, which is \eqref{young}.
The general case follows immediately.
Moreover, this proof shows that the inequality in
\eqref{young} is strict unless $uv\ge0$ and $\abs{u}^p=\abs{v}^q$,
in which case $uv=\abs{u}^p=\abs{v}^q$.

From this one can proceed to prove the H\"older inequality
in the normalized case:
\[
  \Re \int_\Omega uv\,d\mu\le1,\qquad \norm{u}_p=\norm{v}_q=1,
\]
with equality if and only if $uv=\abs{u}^p=\abs{v}^q$ a.e.

Recall that a normed space $X$ is called \emph{uniformly convex}
\cite{clarkson} if for each $\epsilon>0$ there is some $\delta>0$ so
that $x,y\in X$, $\norm{x}\le1$, $\norm{y}\le1$ and
$\norm{x+y}>2-\delta$ imply $\norm{x-y}<\epsilon$.

Perhaps a bit more intuitive is the following equivalent condition,
which we might call \emph{thin slices of the unit ball are small}:
Given $\phi\in X^*$ with $\norm{\phi}=1$, define the $\delta$-\emph{slice}
$S_{\phi,\delta}=\{x\in X\colon\norm{x}\le1\text{ and }\Re \phi(x)>1-\delta\}$.
The ``thin slices'' condition states that
for each $\epsilon>0$ there is some $\delta>0$ so that,
if $\phi\in X^*$ with $\norm{\phi}=1$, then
$\norm{x-y}<\epsilon$ for all $x,y\in S_{\phi,\delta}$.

This condition follows trivially from uniform convexity.
The proof of the converse requires a minor trick:
Given $x,y\in X$ with $\norm{x}\le1$, $\norm{y}\le1$ and
$\norm{x+y}>2-\delta$,
pick $\phi\in X^*$ with $\norm{\phi}=1$ and
$\Re\phi(x+y)>2-\delta$.
Then $\Re\phi(x)=\Re\phi(x+y)-\Re\phi(y)>2-\delta-1=1-\delta$,
and similarly for $y$.
If $\delta$ was chosen according to the ``thin slices''
condition, $\norm{x-y}<\epsilon$ follows.

It should be noted that for the above proof to work,
it may be sufficient to prove the ``thin slices'' condition for \emph{some}
functionals.  In fact, we shall only use functionals on $L^p$ arising
from functions in $L^q$, without using the knowledge that all
functionals on $L^p$ are of this form.

If we let $\delta\to0$ then $S_{\phi,\delta}$ shrinks to the set of unit
vectors satisfying $\Re \phi(x)=1$,
and the ``thin slices'' condition guarantees that there can
be only one such $x$.
In the $L^p$ case, this corresponds to the uniqueness, given $v\in L^q$,
of a unit vector $u\in L^p$ yielding equality in H\"older's inequality.

Since this uniqueness came about due to the sharp inequality in
\eqref{young}, it is reasonable to expect that an improved lower
estimate on the difference between the two sides of \eqref{young}
might yield a version of the ``thin slices'' condition.

This is precisely what the present proof does.
In fact, we shall show that given $\epsilon>0$, a single $\delta$ will
suffice for all $L^p$ spaces (for a fixed $p$).

The requisite inequality turns out to be quite difficult to find in
general, so we concentrate on the case $v=1$ in \eqref{young},
and use the improved inequality to show the ``thin slices''
condition for the functional arising from the constant function $1$ on
$L^p$ over a probability space.
Afterwards, in Lemma \ref{lemma:general},
we show that this is sufficient to cover the general case.

\begin{lemma} \label{lemma:normal}
Given $1<p<\infty$ and $\epsilon>0$,
there exists $\delta>0$ so that,
for every probability space $(\Omega,\nu)$
and every measurable function $z$ on $\Omega$,
$\norm{z}_p\le1$ and $\Re\int_\Omega z\,d\nu>1-\delta$
imply $\norm{z-1}_p<\epsilon$.
\end{lemma}
\begin{proof}
Consider the function
\[
  f(u)=\abs{u}^p-1+p(1-\Re u)
\]
and note that $f(u)>0$ everywhere except for the value $f(1)=0$.
(This is the case $v=1$ in \eqref{young}.)
Further, note that $f(u)$ and $\abs{u-1}^p$ are asymptotically equal
as $\abs{u}\to\infty$.
Thus, given $\epsilon>0$, we can find some $\alpha>1$ so that
\[
  \abs{u-1}^p\le\alpha f(u)
  \text{ whenever }
  \abs{u-1}\ge\epsilon.
\]
Assume that $z$ satisfies the stated conditions, and let
$E=\{\omega\in\Omega\colon\abs{z(\omega)-1}<\epsilon\}$.
Then
\begin{align*}
  \norm{z-1}_p^p
  &=\int_E\abs{z-1}^p\,d\nu+\int_{\Omega\setminus E}\abs{z-1}^p\,d\nu
\\&\le\epsilon^p+\alpha\int_\Omega f(z)\,d\nu
\\&\le\epsilon^p+p\alpha\Bigl(1-\int_\Omega\Re z\,d\nu\Bigr)
\\&<\epsilon^p+p\alpha\delta.
\end{align*}
Thus picking $\delta=\epsilon^p/(p\alpha)$ is sufficient to guarantee
$\norm{z-1}_p<2^{1/p}\epsilon$.
\end{proof}

\begin{lemma} \label{lemma:general}
Given $1<p<\infty$, $p^{-1}+q^{-1}=1$ and $\epsilon>0$,
there exists $\delta>0$ so that
the following holds:
If $u$, $w$ are measurable functions on a measure space $\Omega$
with $\norm{u}_p\le1$ and $\norm{w}_q=1$
and  $\int_\Omega\Re uw\,d\mu>1-\delta$,
then $\norm{u-v}_p<\epsilon$,
where $v$ is the function satisfying $vw=\abs{v}^p=\abs{w}^q$ a.e.
\end{lemma}

\begin{proof}
Let $p$ and $\epsilon$ be given, and choose $\delta$ as in
Lemma \ref{lemma:normal}.

Let $u$, $v$ and $w$ be as stated above.
Since nothing is changed by multiplying $u$, $v$ by a complex function
of absolute value $1$, and dividing $w$ by the same function,
we may assume without loss of generality that $v\ge0$ and $w\ge0$.

Let $z=u/v$ where $v\ne0$ and $z=0$ where $v=0$.
Thus $zv=u$ where $v\ne0$ and $zv=0$ where $v=0$.
Since $(u-zv)zv=0$ we find
$\norm{u}_p^p=\norm{u-zv}_p^p+\norm{zv}_p^p$.
Also
$\Re\int_\Omega zvw\,d\mu=\Re\int_\Omega uw\,d\mu>1-\delta$,
so $\norm{zv}_p>1-\delta$,
and $\norm{u-zv}_p^p<1-(1-\delta)^p$.

Let $\nu$ be the probability measure
\[
  d\nu=vw\,d\mu=v^p\,d\mu=w^q\,d\mu.
\]
We find
\[
  \int_\Omega \abs{z}^p\,d\nu
  =\int_{v\ne0} \abs{u}^p\,d\mu\le1,
  \quad
  \Re\int_\Omega z\,d\nu=\Re\int_\Omega uw\,d\mu>1-\delta.
\]
By Lemma \ref{lemma:normal}, we now get
\[
  \epsilon^p
  >\int_\Omega \abs{z-1}^p\,d\nu
  =\int_\Omega \abs{z-1}^pv^p\,d\mu
  =\int_{v\ne0} \abs{u-v}^p\,d\mu.
\]
On the other hand,
\[
  \int_{v=0} \abs{u-v}^p\,d\mu
  =\int_\Omega \abs{(u-zv}^p\,d\mu
  <1-(1-\delta)^p.
\]
We therefore get $\norm{u-v}_p^p<\epsilon+1-(1-\delta)^p$,
and the proof is complete.
\end{proof}

\begin{theorem}[Clarkson]
$L^p$ is uniformly convex when $1<p<\infty$.
\end{theorem}
\begin{proof}
Consider $x, y\in L^p$ with $\norm{x}_p=\norm{y}_p=1$ and
$\norm{x+y}_p>2-\delta$.
Let $v=(x+y)/\norm{x+y}_p$,
and choose $w\in L^q$ with $vw=\abs{v}^p=\abs{w}^q$.
In particular $\norm{v}_p=\norm{w}_q=1$.
Then
\[
  \int_\Omega (x+y)w\,d\mu=\norm{x+y}_p\int_\Omega vw\,d\mu
  =\norm{x+y}_p>2-\delta.
\]
Since also $\Re\int_\Omega yw\,d\mu\le1$,
this implies $\Re\int_\Omega xw\,d\mu>1-\delta$.
If $\delta$ was chosen according to Lemma \ref{lemma:general},
we get $\norm{x-v}_p<\epsilon$.
Similarly $\norm{y-v}_p<\epsilon$,
and so $\norm{x-y}_p<2\epsilon.$
\end{proof}

One reason for our interest in including the uniform convexity in a
standard functional analysis class is that this implies the
reflexivity of these spaces, by the Milman--Pettis theorem.
(For a remarkably brief proof, see \cite{ringrose}.)
This, in turn, can be used to prove the standard duality theorem for
$L^p$ and $L^q$.
Of course, this requires some comparatively heavy machinery,
but it is machinery that is usually included in such classes anyway.

We finish by outlining the proof.
If $v\in L^p$ then $\norm{v}_p=\max_{\norm{w}_q=1}\Re\int_\Omega vw\,d\mu$.
This is part H\"older's inequality,
and part -- assuming we normalize $v$ --
the choice of $w$ satisfying $vw=\abs{v}^p=\abs{w}^q$.
Thus, with the standard duality,
$L^q$ is isometrically embedded in $(L^p)^*$.
Now assuming that $L^q\ne(L^p)^*$,
an appeal to the Hahn--Banach theorem produces
a nonzero bounded linear functional $f$ on $(L^p)^*$ which
vanishes on $L^q$.
Since $L^p$ is reflexive,
$\phi$ is of the form $\phi(y)=y(u)$ for some $u\in L^p$.
In particular, for each $w\in L^q$ we get
$0=\phi(w)=\int_\Omega uw\,d\mu$.
But then $u=0$, which is a contradiction.


\end{document}